\documentclass[titlepage]{article}
\usepackage{mathrsfs}
\usepackage{pifont}
\usepackage[centertags]{amsmath}
\usepackage{amsfonts}
\usepackage{amssymb}
\usepackage{amsthm}
\usepackage{amscd,amssymb,amsthm}
\theoremstyle{definition}

\theoremstyle{remark}

\numberwithin{equation}{section}

\begin{document}
\title{Criterions of Wiener type for minimally thin sets and  rarefied sets associated with the
stationary Schr\"{o}dinger operator in a cone}
\author{Pinhong Long \ \ Zhiqiang Gao\ \ Guantie Deng\thanks{Corresponding author.
\newline ~Email: hnlphlch10601@sina.com(P. H. Long),
gaozhiq@gmail.com(Z. Q. Gao), \newline denggt@bnu.edu.cn(G. T.Deng).
 ~Supported by SRFDP(No. 20060027023) and NSF of China(No.
10671022 and No. 11101039).}
\\{ School of Mathematical Science, Beijing Normal University,}
\\{ Laboratory of Mathematics and Complex Systems,}
\\{ Beijing~100875, P.R.China} }
\date{February 2012}
\maketitle

\begin{abstract}
In the paper we give some criterions for a-minimally thin sets and
a-rarefied sets associated with the stationary Schr\"{o}dinger
operator at a fixed Martin boundary point or $\infty$ with respect
to a cone. Moreover, we show that a positive superfunction on a cone
behaves regularly outside a-rarefied set. Finally we illustrate the
relation between a-minimally thin set and a-rarefied set in a cone.

\noindent{\textit{AMS 2010 Subject Classification}}: 31B05, 31B25,
31C35.

\noindent{\textit{Key words and phrases}}: Stationary
Schr\"{o}dinger operator; superfunction; minimally thin; rarefied
sets; cone.
\end{abstract}

\maketitle \numberwithin{equation}{section}
\newtheorem{theorem}{\bf Theorem}
\newtheorem{corollary}{\bf Corollary}
\newtheorem{lemma}{\bf Lemma}
\newtheorem{definition}[theorem]{Definition}
\newtheorem{remark}{\bf Remark}
\newpage

\setcounter{section}{0}
\section{Introduction}

\hspace {5mm}In the paper we are concerned with some properties for
the generalized subharmonic function associated with the stationary
Schr\"{o}dinger operator. we will give the characterization for
minimally thin sets and rarefied sets about these generalized
subharmonic function. The study of the minimal thinnes had been
exploited a little and attracted many mathematicians. In [14]
Lelong-Ferrand firstly studied the thinness at boundary points for
the subharmonic functions, then L.Na\"{\i}m in [18] obtained the
equivalent conditions as the criterions for minimally thin sets at a
fixed boundary point with respect to $T_n$. D.H.Armatage and
S.J.Gardiner in [4] systematically reviewed these theoretical
production. Based on the research of Lelong-Ferrand, in [9]
M.Ess\'{e}n and H.L.Jackson knew the criterions for minimally thin
sets at $\infty$ with respect to $T_n$ and introduced the criterions
for rarefied sets at $\infty$ with respect to $T_n$. Since a
rarefied set at $\infty$ with respect to $T_n$ is also minimally
thin set at $\infty$ with respect to $T_n$, M.Ess\'{e}n and
H.L.Jackson in [9] extended Lelong-Ferrand's theorem (from [14]) to
much stronger conclusion. They remarked that by these criterions
there exist some connection between both exceptional sets. Next
I.Miyamoto and H.Yoshida in [17] extended these results(from
M.Ess\'{e}n and H.L.Jackson [9]) from $T_n$ to the conical case. In
view of the above statement, we shall give some criterions of Wiener
type for a-minimally thin sets at a fixed Martin boundary point with
respect to $C_n(\Omega)$ which extend L.Na\"{\i}m's results(see [4
or 16]). Similarly we shall give some criterions of Wiener type for
a-minimally thin sets and a-rarefied sets at $\infty$ with respect
to $C_n(\Omega)$. Moreover, we generalize some theorems (Theorems 3
and 4 from  I.Miyamoto and H.Yoshida[17])to the results associated
with the stationary Schr\"{o}dinger operator which are the crucial
part in our paper. To state our results, we will need some notations
and background materials below.

Let ${\bf R}^{n}(n\geq2)$ be the $n$-dimensional Euclidean space and
$\bf S$ its an open set. The boundary and the closure of $\bf S$ are
denoted by $\partial{\bf S}$ and $\overline{\bf S}$, respectively.
In cartesian coordinate a point $P$ is denoted by $(X,x_n),$ where
$X=(x_1,x_2,\ldots,x_{n-1}).$ Let $|P|$ be the Euclidean norm of $P$
and $|P-Q|$ the Euclidean distance of two points $P$ and $Q$ in
${\bf R}^{n}$. The unit sphere and the upper half unit sphere are
denoted by ${\bf S}^{n-1}$ and ${\bf S}_{+}^{n-1}$, respectively.
For $P\in {\bf R}^{n}$ and $r>0$, let $B(P, r)$ be the open ball of
radius $r$ centered at $P$ in ${\bf R}^{n}$, then
$S_{r}=\partial{B(O, r)}$. Furthermore,we denote by $dS_{r}$ the
$(n-1)$-dimensional volume elements induced by the Euclidean metric
on $S_{r}$.

 We introduce the system of spherical coordinates for $P=(X,y)$ by the following formulas
 $$x_1=r\prod^{n-1}_{j=1}\sin\theta_j~~(n\geq2),~~y=r\cos\theta_1$$
 and if $n\geq3,$
 $$x_{n-k+1}=r\cos\theta_k\prod^{k-1}_{j=1}\sin\theta_j~~(2\leq n\geq n-1),$$
 where $0\leq r<\infty,$ $0\leq \theta_j\leq\pi(1\leq j\leq n-2)$ and
 $-\frac{\pi}{2}\leq \theta_{n-1}\leq\frac{3\pi}{2}.$

 Relative to the system of spherical coordinates, the Laplace operator $\Delta$ may be written
 $$\Delta=\frac{n-1}{r}\frac{\partial}{\partial r}+\frac{\partial^2}{\partial r^2}+\frac{\Delta^*}{r^2},$$
 where the explicit form of the Beltrami operator $\Delta^*$ is given by V.Azarin(see[2]).

Let $D$ be an arbitrary domain in ${\bf R}^{n}$ and $\mathscr{A}_a$
denotes the class of nonnegative radial potentials $a(P)$, i.e.
$0\leq a(P)=a(r)$, $P=(r,\Theta)\in D$, such that $a\in
L_{loc}^{b}(D)$ with some $b> {n}/{2}$ if $n\geq4$ and with $b=2$ if
$n=2$ or $n=3$.

If $a\in \mathscr{A}_a$, then the stationary Schr\"{o}dinger
operator with a potential $a(P)$
$$\mathcal{L}_a=-\Delta+a(P)I\eqno(1.1)$$
can  be extended in the usual way from the space $C_0^{\infty}(D)$
to an essentially self-adjoint operator on $L^{2}(D),$ where
$\Delta$ is the Laplace operator and $I$ the identical operator(see
[21, Chap. 13] ). Then $\mathcal{L}_a$ has a Green $a$-function
$G_{D}^{a}(P,Q)$. Here $G_{D}^{a}(P,Q)$ is positive on $D$ and its
inner normal derivative $\partial G_{D}^{a}(P,Q)/{\partial n_Q}$ is
not negative, where ${\partial}/{\partial n_Q}$ denotes the
differentiation at $Q$ along the inward normal into $D$. We write
this derivative by $PI_{D}^{a}(P,Q)$, which is called the Poisson
$a$-kernel with respect to $D$. Denote by $G_{D}^{0}(P,Q)$ the Green
function of Laplacian. It is well known that, for any potential
$a(P)\geq0$,
$$G_{D}^{a}(P,Q)\leq G_{D}^{0}(P,Q).\eqno(1.2)$$ The inverse
inequality is much more elaborate if $D$ is a bounded domain in
${\bf R}^{n}$. M.Cranston, E.Fabes and Z.Zhao (see [6], the case
$n=2$ is implicitly contained in [7]) have proved
$$G_{D}^{a}(P,Q)\geq M(D)G_{D}^{0}(P,Q),\eqno(1.3)$$
where $D$ is a bounded domain, a constant $M(D)=M(D, a(P))$ is
positive and does not depend on points $P$ and $Q$ in $D$. If $a=0$,
then obviously $M(D)\equiv1$.

 Suppose that a function $u\not\equiv -\infty$ is upper
semi-continuous in $D$. We call $u\in[-\infty,+\infty)$ a
subfunction of the Schr\"{o}dinger operator $\mathcal{L}_a$ if at
each point $P\in D$ with $0<r<r(P)$ the generalized mean-value
inequality
$$u(P)\leq \int_{S(P,r)}u(Q)\frac{\partial G_{B(P,r)}^{a}(P,Q)}{\partial n_Q}d\sigma(Q)\eqno(1.4)$$
is satisfied, where $S(P,r)=\partial{B(P,r)}$, $G_{B(P,r)}^{a}(P,Q)$
is the Green $a$-function of $Sch_a$ in $B(P,r)$ and $d\sigma(Q)$
the surface area element on $S(P,r)$(see [20]).

The class of subfunctions in $D$ is denoted by $SbH(a,D)$. If $-u\in
SbH(a,D)$, then we call $u$ a superfunction and denote the class of
superfunctions by $SpH(a,D)$. If a function $u$ is both subfunction
and superfunction, it is clearly continuous and is called an
$a$-harmonic function associated with the operator $Sch_a$. The
class of $a$-harmonic functions is denoted by $H(a,D)=SbH(a,D)\cap
SpH(a,D)$. In terminology we follow B.Ya.Levin and A.Kheyfits (see
[11], [13] and [15]). For simplicity, a point $(1,\Theta)$ on ${\bf
S}^{n-1}$ and the set $\{\Theta; (1,\Theta)\in \Omega\}$ for a set
$\Omega$ ($\Omega\subset {\bf S}^{n-1}$) are often identified with
$\Theta$ and $\Omega$, respectively. For two sets $\Xi\subset {\bf
R}_+$ and $\Omega\subset {\bf S}^{n-1},$ the set
$\{(r,\Theta)\in{\bf R}^{n}; r\in\Xi,(1,\Theta)\in \Omega\}$ in
${\bf R}^{n}$ is simply denoted by $\Xi\times \Omega.$ In
particular, the half space ${\bf R}_{+}\times {\bf
S}_{+}^{n-1}=\{(X, x_n)\in{\bf R}^{n}; x_n>0\}$ will be denoted by
${\bf T}_n$. By $C_n(\Omega)$ we denote the set ${\bf R}_+\times
\Omega$ in ${\bf R}^{n}$ with the domain $\Omega$ on ${\bf S}^{n-1}$
and call it a cone. We mean the sets $I\times\Omega$ and $I\times
\partial{\Omega}$ with an interval on $\bf R_+$ by $C_n(\Omega;I)$
and $S_n(\Omega;I)$, and $C_n(\Omega)\cap S_{r}$ by $S_n(\Omega;
r)$. By $S_n(\Omega)$ we denote $S_n(\Omega; (0,+\infty))$, which is
$\partial{C_n(\Omega)}-\{O\}.$ From now on, we always assume
$D=C_n(\Omega)$ and write $G_{\Omega}^{a}(P,Q)$ instead of
$G_{C_n(\Omega)}^{a}(P,Q)$.

Let $\Omega$ be a domain on ${\bf S}^{n-1}$ with smooth boundary and
$\lambda$ the least positive eigenvalue for $\Delta^\ast$ on
$\Omega$ (see [22, p. 41])
$$(\Delta^\ast+\lambda)\varphi(\Theta)=0 \hspace {3mm} \textrm{on} \hspace {2mm}
\Omega,\eqno(1.5)$$
$$\hspace {17mm}\varphi(\Theta)=0 \hspace {3mm}
\textrm{on} \hspace {2mm}
\partial{\Omega}.$$
The corresponding eigenfunction is denoted by $\varphi(\Theta)$
satisfying $\int_\Omega\varphi^2(\Theta)dS_1=1$. In order to ensure
the existence of $\lambda$ and a smooth $\varphi(\Theta)$, We put a
rather strong assumption on $\Omega$: if $n\geq3,$ then $\Omega$ is
a $C^{2,\alpha}$-domain $(0<\alpha<1)$ on ${\bf S}^{n-1}$ surrounded
by a finite number of mutually disjoint closed
hypersurfaces~(e.g.~see [12, p. 88-89] for the definition of
$C^{2,\alpha}$-domain).

Solutions of an ordinary differential equation
$$-Q''(r)-\frac{n-1}{r}Q'(r)+\left( \frac{\lambda}{r^2}+a(r)\right )Q(r)=0, \hspace {4mm} for~0<r<\infty \
\eqno{(1.6)}$$ are known (see[24] for more references) that if the
potential $a\in \mathscr{A}_a$. We know the equation (1.6) has a
fundamental system of positive solutions $\{V,W\}$ such that $V$ is
nondecreasing with
$$0\leq V(0+)\leq V(r) \hspace {2mm} \textrm{as} \hspace {2mm} r\rightarrow+\infty\eqno{(1.7)}$$
and $W$ is monotonically decreasing with
$$+\infty=W(0+)>W(r)\searrow0 \hspace {2mm} \textrm{as} \hspace {2mm} r\rightarrow+\infty.\eqno{(1.8)}$$
We remark that both $V(r)\varphi(\Theta)$ and $W(r)\varphi(\Theta)$
are harmonic on $C_n(\Omega)$ and vanish continuously on
$S_n(\Omega)$.

We will also consider the class $\mathscr{B}_a$, consisting of the
potentials $a\in \mathscr{A}_a$ such that there exists the finite
limit $\lim\limits_{r\rightarrow\infty}r^2a(r)=\kappa\in[0,\infty)$,
moreover, $r^{-1}|r^2 a(r)-\kappa|\in L(1,\infty)$. If $a\in
\mathscr{B}_a$, then the (super)subfunctions are continuous (e.g.
see [23]). For simplicity, in the rest of paper we assume that
$a\in\mathscr{B}_a$ and we shall suppress this assumption.

Denote
$$\iota_{\kappa}^{\pm}=\frac{2-n\pm\sqrt{(n-2)^2+4(\kappa+\lambda)}}{2},$$
then the solutions $V(r)$ and $W(r)$ to the equation (1.6)
normalized by $V(1)=W(1)=1$ have the asymptotic(see[12])
$$V(r)\approx r^{\iota_{\kappa}^{+}},~~~~W(r)\approx r^{\iota_{\kappa}^{-}},~~~\textrm{as}~~~r\rightarrow\infty\eqno{(1.9)}$$
and
$$\chi=\iota_{\kappa}^{+}-\iota_{\kappa}^{-}=\sqrt{(n-2)^2+4(\kappa+\lambda)},~~~~\chi'=\omega(V(r),W(r))\mid_{r=1},\eqno{(1.10)}$$
where $\chi'$ is their Wronskian at $r=1$.

\vspace{0.3cm} {\bf Remark 1.} If $a=0$ and $\Omega={\bf
S}_{+}^{n-1}$, then $\iota_{0}^{+}=1$, $\iota_{0}^{-}=1-n$ and
$\varphi(\Theta)=(2n s_n^{-1})^{1/2}cos\theta_1,$ where $s_n$ is the
surface area $2\pi^{n/2}\{\Gamma(n/2)\}^{-1}$ of ${\bf S}^{n-1}$.

R.S.Martin introduced the family with parameter of functions that
was called Martin functions later(see M.Brelot[5] or
R.S.Martin[16]). Now we introduce the Martin function $M^a_{\Omega}$
associated with the stationary Schr\"{o}dinger operator as follows

The function $M^a_{\Omega}$ defined on $C_n(\Omega)\times
C_n(\Omega)-(P_0, P_0)$ by
$$M^{a}_{\Omega}(P,Q)=\frac{G^{a}_{\Omega}(P,Q)}{G^{a}_{\Omega}(P_0,Q)}$$
is called the generalized Martin Kernel of $C_n(\Omega)$ (relative
to $P_0$). If $Q=P_0$, the above quotient is interpreted as 0(for
a=0, we refer to D.H.Armitage and S.Gardiner[4]).

\section{Statements of main results}
First we remark that
$$C_1V(r)W(t)\varphi(\Theta)\varphi(\Phi)\leq G^{a}_{\Omega}(P,Q)\leq C_2V(r)W(t)\varphi(\Theta)\varphi(\Phi)\eqno{(2.1)}$$
or
$$C_1V(t)W(r)\varphi(\Theta)\varphi(\Phi)\leq G^{a}_{\Omega}(P,Q)\leq C_2V(t)W(r)\varphi(\Theta)\varphi(\Phi)\eqno{(2.2)}$$
for any $P=(r, \Theta)\in C_n(\Omega)$ and any $Q=(t, \Phi)\in
C_n(\Omega)$ satisfying $0<\frac{r}{t}\leq \frac{4}{5}$ or
$0<\frac{t}{r}\leq \frac{4}{5},$ where $C_1$ and $C_2$ are two
positive constants(See A.Escassut, W.Tutschke and C.C.Yang[11, Chap.
11] and for $a=0$, see V.S.Azarin [2, Lemma 1], M.Ess\'{e}n and
J.L.Lewis [10, Lemma 2]). It is known that the Martin boundary
$\vartriangle$ of $C_n(\Omega)$ is the set $\partial
C_n(\Omega)\cup\{\infty\}$. When we denote the Martin kernel
associated with the stationary Schr\"{o}dinger operator by
$M^a_{\Omega}(P, Q)(P\in C_n(\Omega), Q\in\partial
C_n(\Omega)\cup\{\infty\})$ with respect to a reference point chosen
suitably, for any $P\in C_n(\Omega),$ we see
$$M^a_{\Omega}(P, \infty)=V(r)\varphi(\Theta),~~M^a_{\Omega}(P, O)=\kappa W(r)\varphi(\Theta),$$
where $O$ denotes the origin of $\bf R^n$ and $\kappa$ is a positive
constant. Let $E$ be a subset  of $C_n(\Omega)$ and $u\geq0$ be a
superfunction on $C_n(\Omega)$. The reduced function of $u$ is
defined as
$$R^E_{u}(P)=\inf\{\upsilon(P): \upsilon\in\Phi^E_u\},$$
where
$$\Phi^E_u=\{\upsilon\in SpH(a, C_n(\Omega)): u\geq0~and~\upsilon\geq u~on~E\}.$$
By $\widehat{R}^E_u$ we denote the regularized reduced function of
$u$ relative to $E$, here
$$\widehat{R}^E_u(P)=\liminf_{P'\rightarrow P}R^E_u(P').$$
Compared with the regularization of a superharmonic function, we
easily know that $\widehat{R}^E_u$ is a superfunction on
$C_n(\Omega)$.

If $E\subseteq C_n(\Omega)$ and $Q\in\vartriangle$, then the Riesz
decomposition and the generalized Martin representation allow us to
express $\widehat{R}^{E}_{M^a_{\Omega}(., Q)}$ uniquely in the form
$G^a_{\Omega}\mu+M^a_{\Omega}\nu$, where $G^a_{\Omega}\mu$ and
$M^a_{\Omega}\nu$ are the generalized Green potential and
generalized Martin representation. We say that $E$ is a-minimally
thin at $Q$ with respect to $C_n(\Omega)$ if $\nu(\{Q\})=0$. at last
We remark that $\vartriangle_0=\{Q\in\vartriangle: C_n(\Omega)~is
~a-minimally~thin~at~Q\}$, where $\vartriangle$ is the Martin
boundary of $C_n(\Omega)$. Next we start to sate our main theorems

\vspace{0.3cm} {\bf Theorem 1.} Let $E\subseteq C_n(\Omega)$ and a
fixed point $Q\in\vartriangle\diagdown\vartriangle_0$. The following
are equivalent

(a) $E$ is a-minimally thin at $Q$;

(b) $\widehat{R}^E_{M^a_{\Omega}(., Q)}\neq M^a_{\Omega}(., Q)$;

(c) $\inf\{\widehat{R}^{E\cap\omega}_{M^a_{\Omega}(., Q)}:
\omega~is~a~Martin~topology ~neighbourhood~of~Q\}=0.$

 If $u$ is a positive superfunction, then we shall write $\mu_u$ for
 the measure appearing in the generalized Martin representation of
 the greatest a-harmonic minorant of $u$.

\vspace{0.3cm} {\bf Theorem 2.} Let $E\subseteq C_n(\Omega)$ and a
fixed point $Q\in\vartriangle\diagdown\vartriangle_0$. Suppose that
$Q$ is a Martin topology limit of $E$. The following are equivalent

(a) $E$ is a-minimally thin at $Q$;

(b) there exists a positive superfunction $u$ such that
$$\liminf_{P\rightarrow Q, P\in E}\frac{u(P)}{M^a_{\Omega}(P, Q)}>\mu_u(\{Q\});$$

(c) there is a a-potential $u$ on $C_n(\Omega)$ such that
$$\frac{u(P)}{M^a_{\Omega}(P, Q)}\rightarrow\infty~~(P\rightarrow Q; P\in E).$$

Let $E$ be a bounded subset of $C_n(\Omega)$. Then
$\widehat{R}^E_{M^a_{\Omega}(., \infty)}(P)$ is bounded on
$C_n(\Omega)$ and hence the greatest a-harmonic minorant of
$\widehat{R}^E_{M^a_{\Omega}(., \infty)}(P)$ is zero. By the Riesz
decomposition theorem that there exists a unique positive measure
$\lambda_E$ on $C_n(\Omega)$ such that
$$\widehat{R}^E_{M^a_{\Omega}(., \infty)}(P)=G^a_{\Omega}\lambda_E(P)\eqno{(2.3)}$$
for any $P\in C_n(\Omega)$ and $\lambda_E$ is concentrated on $B_E$
, where
$$B_E=\{P\in C_n(\Omega): E~is~not~thin~at~P\}.$$
For $a=0$, we see M.Brelot[5] and J.L.Doob[8]. According to the
Fatou's lemma, we easily know the condition (b) in Theorem 3 and
(2.6) below.

\vspace{0.3cm} {\bf Theorem 3.} Let $E\subseteq C_n(\Omega)$ and a
fixed point $Q\in\vartriangle\diagdown\vartriangle_0$. Suppose that
$Q$ is a Martin topology limit point of $E$. The following are
equivalent

(a) $E$ is a-minimally thin at $Q$;

(b) there is a a-potential $G^a_{\Omega}\mu$ such that
$$\liminf_{P\rightarrow Q, P\in E}\frac{G^a_{\Omega}\mu(P)}{G^a_{\Omega}(P_0, P)}>\int M^a_{\Omega}(P, Q)d\mu(P);\eqno{(2.4)}$$

(c) there is a a-potential $G^a_{\Omega}\mu'$ such that $\int
M^a_{\Omega}(P, Q)d\mu'(P)<\infty$ and
$$\frac{G^a_{\Omega}\mu'(P)}{G^a_{\Omega}(P_0, P)}\rightarrow\infty~~(P\rightarrow Q; P\in E).\eqno{(2.5)}$$

\vspace{0.3cm} {\bf Theorem 4.} Let $E\subseteq C_n(\Omega)$,
$Q_0\in C_n(\Omega)$and a fixed point
$Q\in\vartriangle\diagdown\vartriangle_0$. Suppose that $Q$ is a
Martin topology limit point of $E$. Then $E$ is a-minimally thin at
$Q$ if and only if there exists a positive superfunction $u$ such
that
$$\liminf_{P\rightarrow Q, P\in E}\frac{u(P)}{G^a_{\Omega}(Q_0, P)}>\liminf_{P\rightarrow Q}\frac{u(P)}{G^a_{\Omega}(Q_0, P)}.\eqno{(2.6)}$$
The generalized Green energy $\gamma^a_{\Omega}(E)$ of $\lambda_E$
is defined by
$$\gamma^a_{\Omega}(E)=\int_{C_n(\Omega)}(G^a_{\Omega}\lambda_E)d\lambda_E.\eqno{(2.7)}$$
Let $E$ be a subset of $C_n(\Omega)$ and $E_k=E\cap I_k(\Omega)$,
where $I_k(\Omega)=\{P=(r, \Omega)\in C_n(\Omega): 2^k\leq r\leq
2^{k+1}\}.$ The above theorems are concerned with the fixed boundary
points. Next we will consider the cases at infinity.

\vspace{0.3cm} {\bf Theorem 5.} A subset $E$ of $C_n(\Omega)$ is
a-minimally thin at $\infty$ with respect to $C_n(\Omega)$ if and
only if
$$\sum^{\infty}_{k=0}\gamma^a_{\Omega}(E_k) W(2^k)V^{-1}(2^k)<\infty.\eqno{(2.8)}$$

A subset $E$ of $C_n(\Omega)$ is said to be a-rarefied at $\infty$
with respect to $C_n(\Omega)$, if there exists a positive
superfunction $\upsilon(P)$ in $C_n(\Omega)$ such that
$$\inf_{P\in C_n(\Omega)}\frac{\upsilon(P)}{M^a_{\Omega}(P, \infty)}\equiv0$$
and
$$E\subset H_{\upsilon},$$
where
$$H_{\upsilon}=\{P=(r, \Theta)\in C_n(\Omega): \upsilon(P)\geq V(r)\}.$$

\vspace{0.3cm} {\bf Theorem 6.} A subset $E$ of $C_n(\Omega)$ is
a-rarefied at $\infty$ with respect to $C_n(\Omega)$ if and only if
$$\sum^{\infty}_{k=0}W(2^k)\lambda^a_{\Omega}(E_k)<\infty.\eqno{(2.9)}$$

\vspace{0.3cm} {\bf Remark 2.} When $a=0$, Theorems 5 and 6 belong
to I.Miyamoto and H.Yoshida[17]. When $a=0$ and $\Omega=\bf
S^{n-1}_+$, these are exactly the results by H.Aikawa and
M.Ess\'{e}n[3].

Set
$$c(\upsilon, a)=\inf_{P\in C_n(\Omega)}\frac{\upsilon(P)}{M^a_{\Omega}(P, \infty)}\equiv0$$
for a positive superfunction $\upsilon(P)$ on $C_n(\Omega)$. We
immediately know that $c(\upsilon, a)<\infty.$ Actually let $u(P)$
be a subfunction on $C_n(\Omega)$ satisfying
$$\limsup_{P\rightarrow Q, P\in C_n(\Omega)}u(P)\leq0\eqno{(2.10)}$$
for any $Q\in \partial C_n(\Omega)\diagdown\{O\}$ and
$$\sup_{P=(r, \Theta)\in C_n(\Omega)}\frac{u(P)}{V(r)\varphi(\Theta)}=\ell( a)<\infty.\eqno{(2.11)}$$
Then we see $\ell(a)>-\infty$(for $a=0$, see H.Yoshida [25]). If we
apply this to $u=-\upsilon,$ we may obtain $c(\upsilon, a)<\infty.$

\vspace{0.3cm} {\bf Theorem 7.} Let $\upsilon(P)$ be a positive
superfunction on $C_n(\Omega)$. Then there exists a a-rarefied set
$E$ at $\infty$ with respect to $C_n(\Omega)$ such that
$\upsilon(P)V^{-1}(r)$ uniformly converges to $c(\upsilon,
a)\varphi(\Theta)$ on $C_n(\Omega)\diagdown E$ as
$r\rightarrow\infty,$ where $P=(r, \Theta)\in C_n(\Omega).$

From the definition of a-rarefied set, we see that the following
fact: given any a-rarefied set $E$ at $\infty$ with respect to
$C_n(\Omega)$ there exists a positive superfunction $\upsilon(P)$ on
$C_n(\Omega)$ such that $\upsilon(P)V^{-1}(r)\geq1$ on $E$ and
$c(\upsilon, a)=0.$ Hence $\upsilon(P)V^{-1}(r)$ does not converge
to $c(\upsilon, a)\varphi(\Theta)=0$ on $E$ as $r\rightarrow\infty.$

Let $u(P)$ be a subfunction on $C_n(\Omega)$ satisfying (2.10) and
(2.11). Then
$$\upsilon(P)=\ell(a)V(r)\varphi(\Theta)-u(P)~~(P=(r, \Theta)\in C_n(\Omega))$$
is a positive superfunction on $C_n(\Omega)$ such that $c(\upsilon,
a)=0.$ If we apply theorem 7 to this $\upsilon(P)$, then we obtain
the following corollary

\vspace{0.3cm} {\bf Corollary.} Let $u(P)$ be a subfunction on
$C_n(\Omega)$ satisfying (2.10) and (2.11) for $P\in C_n(\Omega)$.
Then there exists a a-rarefied set $E$ at $\infty$ with respect to
$C_n(\Omega)$ such that $\upsilon(P)V^{-1}(r)$ uniformly converges
to $\ell(a)\varphi(\Theta)$ on $C_n(\Omega)\diagdown E$ as
$r\rightarrow\infty,$ where $P=(r, \Theta)\in C_n(\Omega).$

A cone $C_n(\Omega')$ is called a subcone of $C_n(\Omega)$ if
$\overline{\Omega'}\subset\Omega$ where $\overline{\Omega'}$ is the
closure of $\Omega'\subset\bf S^{n-1}$

\vspace{0.3cm} {\bf Theorem 8.} Let $E$ be a subset of
$C_n(\Omega)$. If $E$ is a a-rarefied set at $\infty$ with respect
to $C_n(\Omega)$, then $E$ is a-minimally thin at $\infty$ with
respect to $C_n(\Omega)$. If $E$ is contained in a subcone of
$C_n(\Omega)$ and $E$ is a-minimally thin at $\infty$ with respect
to $C_n(\Omega)$, then $E$ is a-rarefied set at $\infty$ with
respect to $C_n(\Omega)$.

\section{Some Lemmas} \qquad
In our arguments we need the following results

\vspace{0.3cm} {\bf Lemma 1.} Let $E_1, E_2, \cdots, E_m\subseteq
C_n(\Omega)$ and $Q\in\vartriangle.$

(i) If $E_1\subseteq E_2$ and $E_2$ is a-minimally thin at $Q$, then
$E_1$ is a-minimally thin at $Q$;

(ii) If $E_1, E_2, \cdots, E_m$ are a-minimally thin at $Q$, then
$\bigcup^{m}_{k=1}E_k$ is a-minimally thin at $Q$;

(iii) If $E_1$ is a-minimally thin at $Q$, then there is an open
subset $E$ of $C_n(\Omega)$ such that $E_1\subseteq E$ and $E$ is
a-minimally thin at $Q$.

\vspace{0.3mm}{\bf Proof.} Since $\widehat{R}^{E_1}_{M^a_{\Omega}(.,
Q)}\leq\widehat{R}^{E_2}_{M^a_{\Omega}(., Q)}$, we see (i) holds. To
prove (ii) we note that $\widehat{R}^{E_k}_{M^a_{\Omega}(., Q)}$ is
a a-potential for each $k$ and
$$\sum^{m}_{k=1}\widehat{R}^{E_k}_{M^a_{\Omega}(., Q)}\geq M^a_{\Omega}(., Q)~~~~quasi-everywhere~on~\bigcup^{m}_{k=1}E_k,$$
so $\widehat{R}^{\bigcup_{k}E_k}_{M^a_{\Omega}(., Q)}$ is a
a-potential. Finally, To prove (iii), let
$u=\widehat{R}^{E_1}_{M^a_{\Omega}(., Q)}$. Then $u$ is a
a-potential and $u\geq M^a_{\Omega}(., Q)$ on $E_1\setminus F$ for
some polar set $F$. Let $\upsilon$ be a non-zero a-potential such
that $\upsilon=\infty$ on $F$ and let
$$Z=\{P\in C_n(\Omega): u(P)+\upsilon(P)\geq M^a_{\Omega}(P, Q)\}.$$
Then $Z$ is open, $E_1\subseteq Z$ and $R^Z_{M^a_{\Omega}(., Q)}\leq
u+\upsilon$, so $R^Z_{M^a_{\Omega}(., Q)}$ is a a-potential and $Z$
is a-minimally thin at $Q$.

\vspace{0.3cm} {\bf Lemma $2.^{[19]}$}
$$\frac{\partial G^a_{\Omega}(P, Q)}{\partial n_Q}\thickapprox t^{-1}V(t)W(r)\varphi(\Theta)\frac{\partial \varphi(\Phi)}{\partial n_{\Phi}}\eqno{(3.1)}$$
$$(\textrm{resp.}~\frac{\partial G^a_{\Omega}(P, Q)}{\partial n_Q}\thickapprox V(r)t^{-1}W(t)\varphi(\Theta)\frac{\partial \varphi(\Phi)}{\partial n_{\Phi}})\eqno{(3.2)}$$
for any $P=(r,\Theta)\in C_n(\Omega)$ and any $Q=(t,\Phi)\in
S_n(\Omega)$ satisfying $0<\frac{t}{r}\leq \frac{4}{5}$
$(\textrm{resp.}~ 0<\frac{r}{t}\leq \frac{4}{5});$
$$\frac{\partial G^0_{\Omega}(P, Q)}{\partial n_Q}\lesssim \frac{\varphi(\Theta)}{t^{n-1}}\frac{\partial \varphi(\Phi)}{\partial n_{\Phi}}
+\frac{r\varphi(\Theta)}{|P-Q|^{n}}\frac{\partial
\varphi(\Phi)}{\partial n_{\Phi}}\eqno{(3.3)}$$ for any
$P=(r,\Theta)\in C_n(\Omega)$ and any $Q=(t,\Phi)\in S_n(\Omega;
(\frac{4}{5}r,\frac{5}{4}r)).$

\vspace{0.3cm} {\bf Lemma $3.^{[19]}$} Let $\mu$ be a positive
measure on $C_n(\Omega)$ such that there is a sequence of points
$P_i=(r_i, \Theta_i)\in C_n(\Omega)$
$r_i\rightarrow\infty(i\rightarrow\infty)$ satisfying
$$G^a_{\Omega}\mu(P_i)=\int_{C_n(\Omega)}G^a_{\Omega}(P_i, Q)d\mu(t, \Phi)<\infty$$
for $i=1,2,3,\cdots;$ $Q=(t, \Phi)\in C_n(\Omega).$ Then for a
positive number $\ell$
$$\int_{C_n(\Omega; \ell, \infty)}W(t)\varphi(\Phi)d\mu(t, \Phi)<\infty\eqno{(3.4)}$$
and
$$\lim_{R\rightarrow\infty}\frac{W(R)}{V(R)}\int_{C_n(\Omega; 0, R)}V(t)\varphi(\Phi)d\mu(t, \Phi)=0.\eqno{(3.5)}$$

\vspace{0.3cm} {\bf Lemma $4.^{[19]}$} Let $\nu$ be a positive
measure on $S_n(\Omega)$ such that there is a sequence of points
$P_i=(r_i, \Theta_i)\in C_n(\Omega)$
$r_i\rightarrow\infty(i\rightarrow\infty)$ satisfying
$$\int_{S_n(\Omega)}\frac{\partial G^a_{\Omega}(P_i, Q)}{\partial n_Q}d\nu(Q)<\infty$$
for $i=1,2,3,\cdots;$ $Q=(t, \Phi)\in C_n(\Omega).$ Then for a
positive number $\ell$
$$\int_{S_n(\Omega; \ell, \infty)}W(t)t^{-1}\frac{\partial\varphi(\Phi)}{\partial n_{\Phi}}d\nu(t, \Phi)<\infty\eqno{(3.6)}$$
and
$$\lim_{R\rightarrow\infty}\frac{W(R)}{V(R)}\int_{S_n(\Omega; 0, R)}V(t)t^{-1}\frac{\partial\varphi(\Phi)}{\partial n_{\Phi}}d\nu(t, \Phi)=0.\eqno{(3.7)}$$

\vspace{0.3cm} {\bf Lemma 5.} Let $\mu$ be a positive measure on
$C_n(\Omega)$ for which $G^a_{\Omega}\mu(P)$ is defined. Then for
any positive number $A$ the set
$$\{P=(r, \Theta)\in C_n(\Omega): G^a_{\Omega}\mu(P)\geq AV(r)\varphi(\Theta)\}$$
is a-minimally thin at $\infty$ with respect to $C_n(\Omega)$.

\vspace{0.3cm} {\bf Lemma 6.} Let $\upsilon(P)$ be a positive
superfunction on $C_n(\Omega)$ and put
$$c(\upsilon, a)=\inf_{P\in C_n(\Omega)}\frac{\upsilon(P)}{M^a_{\Omega}(P,
\infty)};~~~~c_0(\upsilon, a)=\inf_{P\in
C_n(\Omega)}\frac{\upsilon(P)}{M^a_{\Omega}(P, O)}.\eqno{(3.8)}$$
Then there are a unique positive measure $\mu$ on $C_n(\Omega)$ and
a unique positive measure $\nu$ on $S_n(\Omega)$ such that
$$\upsilon(P)=c(\upsilon, a)M^a_{\Omega}(P, \infty)+c_o(\upsilon, a)M^a_{\Omega}(P, O)
+\int_{C_n(\Omega)}G^a_{\Omega}(P,
Q)d\mu(Q)+\int_{S_n(\Omega)}\frac{\partial G^a_{\Omega}(P,
Q)}{\partial n_Q}d\nu(Q),$$ where $\frac{\partial}{\partial n_Q}$
denotes the differentiation at $Q$ along the inward normal into
$C_n(\Omega)$.

\vspace{0.3mm}{\bf Proof.} By the Riesz decomposition theorem, we
have a unique measure $\mu$ on $C_n(\Omega)$ such that
$$\upsilon(P)=\int_{C_n(\Omega)}G^a_{\Omega}(P, Q)d\mu(Q)+h(P)~~(P\in C_n(\Omega)),\eqno{(3.9)}$$
where $h$ is the greatest a-harmonic minorant of $\upsilon$ on
$C_n(\Omega)$. Furthermore, by the generalized Martin representation
theorem we have with another positive measure $\nu'$ on $\partial
C_n(\Omega)\cup \{\infty\}$
$$h(P)=\int_{\partial C_n(\Omega)\cup \{\infty\}}M^a_{\Omega}(P, Q)d\nu'(Q)$$
$$=M^a_{\Omega}(P, \infty)\nu'(\{\infty\})+M^a_{\Omega}(P, O)\nu'(\{O\})$$
$$+\int_{S_n(\Omega)}M^a_{\Omega}(P, Q)d\nu'(Q)~~(P\in C_n(\Omega)).$$
We know from (3.9) that $\nu'(\{\infty\})=c(\upsilon, a)$ and
$\nu'(\{O\})=c_o(\upsilon, a)$.

Since
$$M^a_{\Omega}(P, Q)=\lim_{P_1\rightarrow Q, P_1\in C_n(\Omega)}\frac{G^a_{\Omega}(P, P_1)}{G^a_{\Omega}(P_0, P_1)}
=\frac{\partial G^a_{\Omega}(P, Q)}{\partial n_Q}/\frac{\partial
G^a_{\Omega}(P_0, Q)}{\partial n_Q},\eqno{(3.10)}$$ where $P_0$ is a
fixed reference point of the Martin kernel, we also obtain
$$h(P)=c(\upsilon, a)M^a_{\Omega}(P, \infty)+c_o(\upsilon, a)M^a_{\Omega}(P, O)$$
$$+\int_{S_n(\Omega)}\frac{\partial G^a_{\Omega}(P, Q)}{\partial n_Q}d\nu(Q)~~(P\in C_n(\Omega))\eqno{(3.11)}$$
by taking
$$d\nu(Q)=\left\{\frac{\partial G^a_{\Omega}(P_0, Q)}{\partial n_Q}\right\}^{-1}d\nu'(Q)~~(Q\in S_n(\Omega)).$$
Hence by (3.9)and (3.11) we get the required.

\vspace{0.3cm} {\bf Lemma 7.} Let $E$ be a bounded subset of
$C_n(\Omega)$ and $u(P)$ be a positive superfunction on
$C_n(\Omega)$ such that $u(P)$ is represented as
$$u(P)=\int_{C_n(\Omega)}G^a_{\Omega}(P, Q)d\mu_u(Q)+\int_{S_n(\Omega)}\frac{\partial G^a_{\Omega}(P,
Q)}{\partial n_Q}d\nu_u(Q)\eqno{(3.12)}$$ with two positive measures
$\mu_u(Q)$ and $\nu_u(Q)$ on $C_n(\Omega)$ and $S_n(\Omega)$
respectively, and satisfies $u(P)\geq1$ for any $P\in E.$ Then
$$\lambda(E)\leq\int_{C_n(\Omega)}V(t)\varphi(\Phi)d\mu_u(t, \Phi)+\int_{S_n(\Omega)}V(t)t^{-1}\varphi(\Phi)d\nu_u(t,\Phi).\eqno{(3.13)}$$
When $u(P)=\widehat{R}^E_1(P)(P\in C_n(\Omega)),$ the equality holds
in (3.13).

\vspace{0.3mm}{\bf Proof.} Since $\lambda_E$ is concentrated on
$B_E$ and $u(P)\geq1$ for any $P\in B_E$, we see that
$$\lambda^a_{\Omega}(E)\leq\int_{C_n(\Omega)}d\lambda_E(P)\leq\int_{C_n(\Omega)}u(P)d\lambda_E(P)$$
$$\leq\int_{C_n(\Omega)}\widehat{R}^E_{M^a_{\Omega}(., \infty)}d\mu_u(Q)$$
$$+\int_{S_n(\Omega)}\left(\int_{C_n(\Omega)}\frac{\partial G^a_{\Omega}(P, Q)}{\partial n_Q}d\lambda_E(P)\right)d\nu_u(Q).\eqno{(3.14)}$$
In addition we have
$$\widehat{R}^E_{M^a_{\Omega}(., \infty)}(Q)\leq M^a_{\Omega}(Q, \infty)=V(t)\varphi(\Phi)~~(Q=(t, \Phi)\in C_n(\Omega)).\eqno{(3.15)}$$
Since
$$\int_{C_n(\Omega)}\frac{\partial G^a_{\Omega}(P, Q)}{\partial n_Q}d\lambda_E(P)
\leq\liminf_{\rho\rightarrow0}\frac{1}{\rho}\int_{C_n(\Omega)}
G^a_{\Omega}(P, P_{\rho})d\lambda_E(P)$$ for any $Q\in S_n(\Omega)$
$(P_{\rho}=(r_{\rho}, \Theta_{\rho})=Q+\rho n_Q\in C_n(\Omega)),$
$n_Q$ is the inward normal unit vector at $Q$ and
$$\int_{C_n(\Omega)} G^a_{\Omega}(P, P_{\rho})d\lambda_E(P)=\widehat{R}^E_{M^a_{\Omega}(., \infty)}(P_{\rho}) \leq
M^a_{\Omega}(P_{\rho},\infty)=V(r_{\rho})\varphi(\Phi_{\rho}),$$ we
have
$$\int_{C_n(\Omega)}\frac{\partial G^a_{\Omega}(P, Q)}{\partial n_Q}d\lambda_E(P)
\leq V(t)t^{-1}\frac{\partial \varphi(\Phi)}{\partial
n_Q}\eqno{(3.16)}$$ for any $Q=(t, \Phi)\in S_n(\Omega)$. Thus we
obtain (3.13). Because $\widehat{R}^E_1(P)$ is bounded on
$C_n(\Omega)$, $u(P)$ has the expression (3.12) by the Lemma 6 when
$u(P)=\widehat{R}^E_1(P)$. Then we easily have the equalities only
in (3.14) because $\widehat{R}^E_1(P)=1$ for any $P\in
B_E$(J.L.Doob[8 p.169]). Hence we claim if
$$\mu_u(\{P\in C_n(\Omega): \widehat{R}^E_{M^a_{\Omega}(., \infty)}(P)<M^a_{\Omega}(P,\infty)\})=0\eqno{(3.17)}$$
and
$$\nu_u(\{Q=(t, \Phi)\in S_n(\Omega): \int_{C_n(\Omega)}\frac{\partial G^a_{\Omega}(P, Q)}{\partial n_Q}d\lambda_E(P)
\leq V(t)t^{-1}\frac{\partial \varphi(\Phi)}{\partial
n_Q}\})=0,\eqno{(3.18)}$$ then the equality in (3.13) holds.

To see (3.17) we remark that
$$\{P\in C_n(\Omega): \widehat{R}^E_{M^a_{\Omega}(., \infty)}(P)<M^a_{\Omega}(P,\infty)\}\subset C_n(\Omega)\diagdown B_E$$
and
$$\mu_u(C_n(\Omega)\diagdown B_E)=0.$$
To prove (3.18) we set
$$B'_E=\{Q\in S_n(\Omega): E~is~not~a-minimally~thin~at~Q\}\eqno{(3.19)}$$
and
$$e=\{P\in E: \widehat{R}^E_{M^a_{\Omega}(., \infty)}(P)<M^a_{\Omega}(P,\infty)\}.\eqno{(3.20)}$$
Then $e$ is a apolar set and hence for any $Q\in S_n(\Omega)$
$$\widehat{R}^E_{M^a_{\Omega}(., Q)}=\widehat{R}^{E\diagdown e}_{M^a_{\Omega}(., Q)}.$$
Consequently, for any $Q\in B'_E$ $E\diagdown e$ is not also
a-minimally thin at $Q$ and
$$\int_{C_n(\Omega)}M^a_{\Omega}(P, Q)d\eta(P)=\lim_{P'\rightarrow Q, P'\in E\diagdown e}
\int_{C_n(\Omega)}M^a_{\Omega}(P, P')d\eta(P)\eqno{(3.21)}$$ for any
positive measure $\eta$ on $C_n(\Omega)$, where
$$M^a_{\Omega}(P, P')=\frac{G^a_{\Omega}(P, P')}{G^a_{\Omega}(P_0, P')}~~(P\in C_n(\Omega), P'\in C_n(\Omega)).$$
Now we take $\eta=\mu_E$ in (3.21). Since
$$\lim_{P\rightarrow Q, P\in C_n(\Omega)}\frac{M^a_{\Omega}(P, \infty)}{G^a_{\Omega}(P_0, P)}=
V(t)t^{-1}\frac{\partial \varphi(\Phi)}{\partial n_{\Phi}} \left\{
\frac{\partial G^a_{\Omega}(P_0, Q)}{\partial
N_Q}\right\}^{-1}~~(Q=(t, \Phi)\in S_n(\Omega)),$$ we obtain from
(3.10)
$$\int_{C_n(\Omega)}\frac{\partial G^a_{\Omega}(P, Q)}{\partial n_{\Phi}}d\mu_E(P)=
V(t)t^{-1}\frac{\partial \varphi(\Phi)}{\partial
n_{\Phi}}\liminf_{P'\rightarrow Q, P'\in E\diagdown
e}\int_{C_n(\Omega)}\frac{G^a_{\Omega}(P, P')}{M^a_{\Omega}(P',
\infty)}d\mu_E(P)\eqno{(3.22)}$$ for any $Q\in (t, \Phi)\in B'_E$.
Since
$$\int_{C_n(\Omega)}\frac{G^a_{\Omega}(P, P')}{M^a_{\Omega}(P', \infty)}d\mu_E(P)=
\frac{1}{M^a_{\Omega}(P', \infty)}\widehat{R}^E_{M^a_{\Omega}(.,
\infty)}(P')=1$$ for any $P'\in E\diagdown e$, we have
$$\int_{C_n(\Omega)}\frac{\partial G^a_{\Omega}(P, Q)}{\partial n_{\Phi}}d\mu_E(P)=
V(t)t^{-1}\frac{\partial \varphi(\Phi)}{\partial n_{\Phi}}$$ for any
$Q=(t, \Phi)\in B'_E$, which shows
$$\{Q=(t, \Phi)\in S_n(\Omega): \int_{C_n(\Omega)}\frac{\partial G^a_{\Omega}(P, Q)}{\partial
n_{\Phi}}d\mu_E(P)<V(t)t^{-1}\frac{\partial \varphi(\Phi)}{\partial
n_{\Phi}}\}\subset S_n(\Omega)\diagdown B'_E.\eqno{(3.23)}$$ Let $h$
be the greatest a-harmonic minorant of $u(P)=\widehat{R}^E_1(P)$ and
$\nu'_u$ be the generalized Martin representing measure of $h$, we
can prove that
$$\widehat{R}^E_h(P)=h\eqno{(3.24)}$$
on $C_n(\Omega)$, then $\nu'_u(S_n(\Omega)\diagdown B'_E)=0$. Since
$$d\nu'_u(Q)=\frac{\partial G^a_{\Omega}(P_0, Q)}{\partial n_Q}d\nu_u(Q)~~(Q\in S_n(\Omega))$$
from (3.10), we also have $\nu_u(S_n(\Omega)\diagdown B'_E)=0,$
which gives (3.18) from (3.23).

To prove (3.24), we set $u^*=\widehat{R}^E_1(P)-h$. Then
$$u^*+h=\widehat{R}^E_1=\widehat{R}^E_{u^*+h}\leq\widehat{R}^E_{u^*}+\widehat{R}^E_{h}$$
and hence
$$\widehat{R}^E_{h}-h\geq u^*-\widehat{R}^E_{u^*}\geq0$$
from which (3.24) follows.

\vspace{0.3cm} {\bf Lemma 8 (The generalized Martin
representation).} If $u$ is a positive a-harmonic function on
$C_n(\Omega)$, then there exists a measure $\mu_u$ on
$\vartriangle$, uniquely determined by $u$, such that
$\mu_u(\vartriangle_0)=0$ and
$$u(P)=\int_{\vartriangle}M^a_{\Omega}(P, Q)d\mu_u(Q)~~~~(P\in C_n(\Omega)),$$
where $\vartriangle_0$ is the same as introduction.

\vspace{0.3cm} {\bf Remark 3.} Following the same method of
D.H.Armitage and S.J.Gardiner[4] for Martin representation we may
easily prove Lemma 8.

\section{Proof of the Theorems}
\vspace{0.3mm}{\bf Proof of the Theorem 1.} First we assume that (b)
holds and let $u=\widehat{R}^{E}_{M^a_{\Omega}(., Q)}$. Since
$M^a_{\Omega}(., Q)$ is minimal, the Riesz decomposition of $u$ is
of the form $\upsilon+\ell M^a_{\Omega}(., Q)$, where $\upsilon$ is
a potential with the stationary Schr\"{o}dinger  operator on
$C_n(\Omega)$ and $0<\ell<1$. Since $u=M^a_{\Omega}(., Q)$
quasi-everywhere on $E$ and $\widehat{R}^{E}_{\upsilon}+\ell
u=\upsilon+\ell M^a_{\Omega}(., Q)=u$ quasi-everywhere on $E$,
$$\widehat{R}^{E}_{M^a_{\Omega}(., Q)}=\widehat{R}^{E}_{u}\leq\widehat{R}^{E}_{\upsilon}+\ell u
\leq\upsilon+\ell M^a_{\Omega}(.,
Q)=\widehat{R}^{E}_{M^a_{\Omega}(., Q)}.\eqno{(4.1)}$$ Hence $\ell
(M^a_{\Omega}(., Q)-u)\equiv0$, so $\ell=0$ by the hypothesis and
(a) holds.

Next we assume (a) holds and let ${\omega_m}$ be a decreasing
sequence of compact neighborhoods of $Q$ in the Martin topology such
that $\bigcap_{m}\omega_m=\{Q\}$. Then $\widehat{R}^{E\bigcap
\omega_m}_{M^a_{\Omega}(., Q)}$ is a-harmonic on
$C_n(\Omega)\diagdown\omega_m$ and the decreasing sequence
$\{\widehat{R}^{E\bigcap \omega_m}_{M^a_{\Omega}(., Q)}\}$ has a
limit $h$ which is a-harmonic on $C_n(\Omega)$. Since $h$ is
majorized by the potential $\widehat{R}^{E}_{M^a_{\Omega}(., Q)}$,
it follows that $h\equiv0$ and (c) follows.

Finally we assume (c) holds, then there is a Martin topology
neighborhood $\omega$ of $Q$ such that $\widehat{R}^{E\bigcap
\omega}_{M^a_{\Omega}(., Q)}\neq M^a_{\Omega}(., Q)$. Since (b)
implies (a), the set $E\bigcap \omega$ is a-minimally thin at $Q$
and so $\widehat{R}^{E\bigcap \omega}_{M^a_{\Omega}(., Q)}$ is a
a-potential. Then $\widehat{R}^{E}_{M^a_{\Omega}(., Q)}$ is a
a-potential and we yields (b).

\vspace{0.3mm}{\bf Proof of the Theorem 2.} Obviously we see that
(c) implies (b). If (b) holds, then there exist $\ell>\mu_u(\{Q\})$
and a Martin topology neighborhood $\omega$ of $Q$ such that $u\geq
\ell M^a_{\Omega}(., Q)$ on $E\bigcap \omega$. If
$\widehat{R}^{E\bigcap\omega}_{M^a_{\Omega}(., Q)}=M^a_{\Omega}(.,
Q)$, then $u\geq \widehat{R}^{E\bigcap\omega}_u\geq\ell
M^a_{\Omega}(., Q)$ and this yields contradictory conclusion that
$\mu_u=\ell\delta_Q+\mu_{u-\ell M^a_{\Omega}(.,
Q)}>\mu_u(\{Q\})\delta_Q$, where $\delta_Q$ is the unit measure with
support $\{Q\}$. Hence
$\widehat{R}^{E\bigcap\omega}_{M^a_{\Omega}(., Q)}\neq
M^a_{\Omega}(., Q)$. Thus $E\bigcap \omega$ is a-minimally thin at
$Q$ and so (a) holds.

Finally we assume (a) holds. By Lemma 1 there is an open subset $U$
of $C_n(\Omega)$ such that $E\subseteq U$ and $U$ is a-minimally
thin at $Q$. By Theorem 1 there is a sequence $\{\omega_m\}$ of
Martin topology open neighborhoods of $Q$ such that
$\widehat{R}^{E\bigcap\omega_m}_{M^a_{\Omega}(., Q)}(P_0)<2^{-m}$.
The function
$u_1=\sum_{n}\widehat{R}^{U\bigcap\omega_m}_{M^a_{\Omega}(., Q)}$,
being a sum of a-potentials, is a a-potential since
$u_1(P_0)<\infty$. Further, since
$\widehat{R}^{E\bigcap\omega_m}_{M^a_{\Omega}(., Q)}=M^a_{\Omega}(.,
Q)$ on the open set $E\bigcap \omega_m$,
$$\frac{u_1(P)}{M^a_{\Omega}(P, Q)}\rightarrow\infty~~(P\rightarrow Q; P\in U)$$
and so (c) holds.

\vspace{0.3mm}{\bf Proof of the Theorem 3.} Clearly (c) implies (b).
To prove that (b) implies (a), we suppose that (b) holds and choose
A such that
$$\liminf_{P\rightarrow Q, P\in E}\frac{G^a_{\Omega}\mu(P)}{G^a_{\Omega}(P_0, P)}>A>\int M^a_{\Omega}(., Q)d\mu.$$
Then $G^a_{\Omega}\mu>AG^a_{\Omega}(P_0, .)$ on $E\cap\omega$ for
some Martin topology neighborhood $\omega$ of $Q$. If $\nu$ denotes
the swept measure of $\delta_{P_0}$ onto $E\cap\omega$, where
$\delta_{P_0}$ is the unit measure with support $\{P_0\}$, then it
follows that
$$G^a_{\Omega}\mu\geq A\widehat{R}^{E\cap\omega}_{G^a_{\Omega}(P_0, .)}=AG^a_{\Omega}\nu$$
on $C_n(\Omega)$. Let $\{K_n\}$ be a sequence of compact subsets of
$C_n(\Omega)$ such that $\bigcup_nK_n=C_n(\Omega)$ and let
$G^a_{\Omega}\mu_n$ denotes the a-potential
$\widehat{R}^{K_n}_{M^a_{\Omega}(., Q)}$. Then
$$\int\widehat{R}^{K_n}_{M^a_{\Omega}(., Q)}d\nu=\int G^a_{\Omega}\nu
d\mu_n \leq A^{-1}\int G^a_{\Omega}\mu
d\mu_n=A^{-1}\int\widehat{R}^{K_n}_{M^a_{\Omega}(., Q)}d\mu.$$
Letting $n\rightarrow\infty$, we see from our choice of $A$ that
$$\widehat{R}^{E\cap\omega}_{M^a_{\Omega}(., Q)}(P_0)=\int M^a_{\Omega}(., Q)d\nu
\leq A^{-1}\int M^a_{\Omega}(., Q) d\mu<1=M^a_{\Omega}(P_0, Q),$$
then $E\cap\omega$ is a-minimally thin at $Q$ by Theorem 1 and so
(a) holds.

Next we suppose that (a) holds. By Lemma 1 there is an open subset
$U$ of $C_n(\Omega)$ such that $E\subseteq U$ and $U$ is a-minimally
thin at $Q$. By Theorem 1 there is a sequence $\{\omega_n\}$ of
Martin topology open neighborhoods of $Q$ such that
$$\sum_{n}\widehat{R}^{U\cap\omega_n}_{M^a_{\Omega}(., Q)}(P_0)<\infty.$$
Let $\mu'=\sum_{n}\nu_n$, where $\nu_n$ is swept measure of
$\delta_{P_0}$ onto $U\cap\omega_n$. Then
$$\int M^a_{\Omega}(P, Q)d\mu'(P)=\sum_{n}\int M^a_{\Omega}(P, Q)d\nu_n(P)=\sum_{n}
\widehat{R}^{U\cap\omega_n}_{M^a_{\Omega}(., Q)}(P_0)<\infty$$ and
(2.5) holds since
$$G^a_{\Omega}\nu_n=\widehat{R}^{U\cap\omega_n}_{G^a_{\Omega}(P_0, .)}=G^a_{\Omega}(P_0, .)$$
on the open set $U\cap\omega_n$, so (c) holds.

\vspace{0.3mm}{\bf Proof of the Theorem 4.} We observe that (2.6) is
dependent of the choice of $Q_0$ since we may multiply across by
$M^a_{\Omega}(Q_0, Q)$. Thus we assume that $Q_0=P_0$, and we claim
that for any a-potential $G^a_{\Omega}\mu$
$$\liminf_{P\rightarrow Q}\frac{G^a_{\Omega}\mu(P)}{G^a_{\Omega}(P_0, P)}=\int M^a_{\Omega}(P, Q)d\mu(P).\eqno{(4.2)}$$
According to Fatou's Lemma, we yield
$$\liminf_{P\rightarrow Q}\frac{G^a_{\Omega}\mu(P)}{G^a_{\Omega}(P_0, P)}\geq\int M^a_{\Omega}(P, Q)d\mu(P).$$
In addition we know that
$$\liminf_{P\rightarrow Q}\frac{G^a_{\Omega}\mu(P)}{G^a_{\Omega}(P_0, P)}<\int M^a_{\Omega}(P, Q)d\mu(P)$$
from Theorem 3 and the fact that $C_n(\Omega)$ is not a-minimally
thin at $Q$. Thus the claim holds.

If $E$ is a-minimally thin at $Q$, we see from (4.2) and the
condition (b) of Theorem 3 that (2.6) holds for some a-potential
$u$. Conversely, if (2.6) holds, then we choose $A$ such that
$$\liminf_{P\rightarrow Q, P\in E}\frac{u(P)}{G^a_{\Omega}(P_0, P)}>A>\liminf_{P\rightarrow Q}\frac{u(P)}{G^a_{\Omega}(P_0, P)}$$
and define $G^a_{\Omega}\mu$ to be the a-potential $\min\{u,
AG^a_{\Omega}(P_0, .)\}$. Then by (4.2)
$$\liminf_{P\rightarrow Q, P\in E}\frac{G^a_{\Omega}\mu(P)}{G^a_{\Omega}(P_0, P)}=A>\liminf_{P\rightarrow Q}\frac{G^a_{\Omega}\mu(P)}{G^a_{\Omega}(P_0, P)}=\int M^a_{\Omega}(P, Q)d\mu(P),$$
and it follows from Theorem 3 that $E$ is a-minimally thin at $Q$.

\vspace{0.3mm}{\bf Proof of the Theorem 5.} We apply the Riesz
decomposition theorem to the superfunction
$\widehat{R}^E_{M^a_{\Omega}(., \infty)}$ on $C_n(\Omega)$, then we
have a positive measure $\mu$ on $C_n(\Omega)$ satisfying
$$G^a_{\Omega}\mu(P)<\infty$$
for any $P\in C_n(\Omega)$ and a non-negative greatest a-harmonic
minorant $H$ of $\widehat{R}^E_{M^a_{\Omega}(., \infty)}$ such that
$$\widehat{R}^E_{M^a_{\Omega}(., \infty)}=G^a_{\Omega}\mu(P)+H.\eqno{(4.3)}$$
We remark that $M^a_{\Omega}(., \infty)$$(P\in C_n(\Omega))$ is a
minimal function at $\infty$. If $E$ is a-minimally thin at $\infty$
with respect to $C_n(\Omega)$, then $\widehat{R}^E_{M^a_{\Omega}(.,
\infty)}$ is a a-potential and hence $H\equiv0$ on $C_n(\Omega)$.
Since
$$\widehat{R}^E_{M^a_{\Omega}(., \infty)}(P)=M^a_{\Omega}(P,
\infty)\eqno{(4.4)}$$ for any $P\in B_E$, we see from (4.3)
$$G^a_{\Omega}\mu(P)=M^a_{\Omega}(P,\infty)\eqno{(4.5)}$$
for any $P\in B_E$. We take a sufficiently large $R$ from Lemma 3
such that
$$C_2\frac{W(R)}{V(R)}\int_{C_n(\Omega;
0, R)}V(t)\varphi(\Phi)d\mu(t, \Phi)<\frac{1}{4}.$$ Then from (2.1
or 2.2)
$$\int_{C_n(\Omega; 0, R)}G^a_{\Omega}(P, Q)d\mu(Q)<\frac{1}{4}M^a_{\Omega}(P,\infty)\eqno{(4.6)}$$
for any $P=(r, \Theta)\in C_n(\Omega)$ and $r\geq\frac{5}{4}r$, and
hence from (4.5)
$$\int_{C_n(\Omega; R, \infty)}G^a_{\Omega}(P, Q)d\mu(Q)\geq\frac{3}{4}M^a_{\Omega}(P,\infty)\eqno{(4.7)}$$
for any $P=(r, \Theta)\in B_E$ and $r\geq\frac{5}{4}r$. Now we
divide $G^a_{\Omega}\mu$ into three parts as follows
$$G^a_{\Omega}\mu(P)=A^{(k)}_1(P)+A^{(k)}_2(P)+A^{(k)}_3(P)~~(P=(r, \Theta)\in C_n(\Omega)),\eqno{(4.8)}$$
where
$$A^{(k)}_1(P)=\int_{C_n(\Omega; 2^{k-1}, 2^{k+2})}G^a_{\Omega}(P, Q)d\mu(Q);$$
$$A^{(k)}_2(P)=\int_{C_n(\Omega; 0, 2^{k-1})}G^a_{\Omega}(P, Q)d\mu(Q)$$
and
$$A^{(k)}_3(P)=\int_{C_n(\Omega; 2^{k+2}, \infty)}G^a_{\Omega}(P, Q)d\mu(Q).$$
Then we shall show that there exists an integer $N$ such that
$$B_E\cap\overline{I_k(\Omega)}\subset\{P=(r, \Theta)\in C_n(\Omega): A^{(k)}_1(P)\geq\frac{1}{4}V(r)\varphi(\Theta)\}~~(k\geq N).\eqno{(4.9)}$$
When we choose a sufficiently large integer $N_1$ by Lemma 3 such
that
$$\frac{W(2^{k})}{V(2^{k})}\int_{C_n(\Omega; 0, 2^{k})}V(t)\varphi(\Phi)d\mu(t, \Phi)<\frac{1}{4C_2}~~(k\geq N_1)$$
and
$$\int_{C_n(\Omega; 2^{k+2}, \infty)}W(t)\varphi(\Phi)d\mu(t, \Phi)<\frac{1}{4C_2}~~(k\geq N_1)$$
for any $P=(r, \Theta)\in\overline{I_k(\Omega)}\cap C_n(\Omega),$ we
have from (2.1 or 2.2) that
$$A^{(k)}_2(P)\leq\frac{1}{4}V(r)\varphi(\Theta)\}~~(k\geq N_1)\eqno{(4.10)}$$
and
$$A^{(k)}_3(P)\leq\frac{1}{4}V(r)\varphi(\Theta)\}~~(k\geq N_1).\eqno{(4.11)}$$
Put
$$N=\max\{N_1, [\frac{\log R}{\log2}]+2\}$$
For any $P=(r, \Theta)\in B_E\cap\overline{I_k(\Omega)}(k\geq N)$.
we have from (4.7), (4.8), (4.10) and (4.11) that
$$A^{(k)}_1(P)\geq\int_{C_n(\Omega; R, \infty)}G^a_{\Omega}(P, Q)d\mu(Q)-A^{(k)}_2(P)-A^{(k)}_3(P)\geq\frac{1}{4}V(r)\varphi(\Theta),$$
which shows (4.9).

Since the measure $\lambda_{E_k}$ is concentrated on $B_{E_k}$ and
$B_{E_k}\subset B_E\cap\overline{I_k(\Omega)}$, finally we obtain by
(4.9) that
$$\gamma^a_{\Omega}(E_k)=\int_{C_n(\Omega)}(G^a_{\Omega}\lambda_{E_k})d\lambda_{E_k}(P)\leq\int_{B_{E_k}}V(t)\varphi(\Phi)d\lambda_{E_k}(r, \Theta)
\leq4\int_{B_{E_k}}A^{(k)}_1(P)d\lambda_{E_k}(P)$$
$$\leq4\int_{C_n(\Omega; 2^{k-1}, 2^{k+2})}\left\{\int_{C_n(\Omega)}G^a_{\Omega}(P, Q)d\lambda_{E_k}(P)\right\}d\mu(Q)$$
$$\leq4\int_{C_n(\Omega; 2^{k-1}, 2^{k+2})}V(t)\varphi(\Phi)d\mu(t, \Phi)~~(k\geq N)$$
and hence
$$\sum^{\infty}_{k=N}\gamma^a_{\Omega}(E_k) W(2^k)V^{-1}(2^k)\lesssim\sum^{\infty}_{k=N}
\int_{C_n(\Omega; 2^{k-1}, 2^{k+2})}W(t)\varphi(\Phi)d\mu(t,\Phi)$$
$$=\int_{C_n(\Omega; 2^{N-1},\infty)}W(t)\varphi(\Phi)d\mu(t,\Phi)<\infty$$ from
Lemma 3, (1.7) and Lemma C.1 in [8 or 10], which gives (2.8).

Next we shall prove the sufficient. Since
$$\widehat{R}^{E_k}_{M^a_{\Omega}(., \infty)}(Q)=M^a_{\Omega}(Q, \infty)$$ for any $Q\in
B_{E_k}$ as in (4.4), we have
$$\gamma^a_{\Omega}(E_k)=\int_{B_{E_k}}M^a_{\Omega}(Q, \infty)d\lambda_{E_k}(Q)\geq V(2^k)\int_{B_{E_k}}\varphi(\Phi)d\lambda_{E_k}(t, \Phi)
~~(Q=(t, \Phi)\in C_n(\Omega))$$ and hence from (2.1 or 2.2), (1.7)
and (1.8)
$$\widehat{R}^{E_k}_{M^a_{\Omega}(., \infty)}(P)\leq C_2V(r)\varphi(\Theta)\int_{B_{E_k}}W(t)\varphi(\Phi)d\lambda_{E_k}(t, \Phi)$$
$$\leq C_2V(r)\varphi(\Theta)V^{-1}(2^k)W(2^k)\gamma^a_{\Omega}(E_k)\eqno{(4.12)}$$
for any $P=(r, \Theta)\in C_n(\Omega)$ and any integer $k$
satisfying $2^k\geq \frac{5}{4}r$. If we define a measure $\mu$ on
$C_n(\Omega)$ by
$$
d\mu(Q)=\left\{\begin{array}{ll}
 \sum^{\infty}_{k=0}d\lambda_{E_k}(Q) & ~~~~~~~~   \ \ \ (Q\in C_n(\Omega; [1, \infty)),  \\
 0
  &~~~~~~~~
 \ \ \ (Q\in C_n(\Omega; (0, 1)),
 \end{array}\right.
$$
then from (2.8) and (4.12)
$$G^a_{\Omega}\mu(P)=\int_{C_n(\Omega)}G^a_{\Omega}(P, Q)d\mu(Q)=\sum^{\infty}_{k=0}
\widehat{R}^{E_k}_{M^a_{\Omega}(.,\infty)}(P)$$ is a finite-valued
superfunction on $C_n(\Omega)$ and
$$G^a_{\Omega}\mu(P)\geq\int_{C_n(\Omega)}G^a_{\Omega}(P, Q)d\lambda_{E_k}(Q)
=\widehat{R}^{E_k}_{M^a_{\Omega}(.,\infty)}(P)=V(r)\varphi(\Theta)$$
for any $P=(r, \Theta)\in B_{E_k}$, and from (2.1 or 2.2)
$$G^a_{\Omega}\mu(P)\geq C'V(r)\varphi(\Theta)\eqno{(4.13)}$$
for any $P=(r, \Theta)\in C_n(\Omega; 0, 1)$, where
$$C'=C_1\int_{C_n(\Omega; \frac{5}{4}, \infty)}W(t)\varphi(\Phi)d\mu(t,\Phi).\eqno{(4.14)}$$
If we set
$$E'=\bigcup^{\infty}_{k=0} B_{E_k}, E_1=E\cap
C_n(\Omega; 0, 1)~and~C=\min(C', 1),$$ then
$$E'\subset\{P=(r, \Theta)\in C_n(\Omega); G^a_{\Omega}\mu(P)\geq CV(r)\varphi(\Theta)\}.$$
Hence by Lemma 5 $E'$ is a-minimally thin at $\infty$ with respect
to $C_n(\Omega)$, namely, there is a point $P'\in C_n(\Omega)$ such
that
$$\widehat{R}^{E'}_{M^a_{\Omega}(., \infty)}(P')\neq M^a_{\Omega}(P', \infty).$$
Since $E'$ is equal  to $E$ except a polar set, we know that
$$\widehat{R}^{E'}_{M^a_{\Omega}(., \infty)}(P)=\widehat{R}^{E}_{M^a_{\Omega}(., \infty)}(P)$$
for any $P\in C_n(\Omega)$ and hence
$$\widehat{R}^{E}_{M^a_{\Omega}(., \infty)}(P')\neq M^a_{\Omega}(P', \infty).$$
This shows that $E$ is a-minimally thin at $\infty$ with respect to
$C_n(\Omega)$.

\vspace{0.3mm}{\bf Proof of the Theorem 6.} Let a subset $E$ of
$C_n(\Omega)$ be a-rarefied set at $\infty$ with respect to
$C_n(\Omega)$, then there exists a positive superfunction
$\upsilon(P)$ on $C_n(\Omega)$ such that $C(\upsilon, a)\equiv0$ and
$$E\subset H_{\upsilon}.\eqno{(4.15)}$$
By Lemma 6 we can find two positive measure $\mu$ on $C_n(\Omega)$
and $\nu$ on $S_n(\Omega)$ such that
$$\upsilon(P)=c_o(\upsilon, a)M^a_{\Omega}(P, O)+\int_{C_n(\Omega)}G^a_{\Omega}(P,
Q)d\mu(Q)+\int_{S_n(\Omega)}\frac{\partial
G^a_{\Omega}(P,Q)}{\partial n_Q}d\nu(Q)~~(P\in C_n(\Omega)).$$ Now
we set
$$\upsilon(P)=c_o(\upsilon, a)M^a_{\Omega}(P, O)+B^{(k)}_1(P)+B^{(k)}_2(P)+B^{(k)}_3(P),\eqno{(4.16)}$$
where
$$B^{(k)}_1(P)=\int_{C_n(\Omega; 0, 2^{k-1})}G^a_{\Omega}(P, Q)d\mu(Q)+\int_{S_n(\Omega; 0, 2^{k-1})}\frac{\partial
G^a_{\Omega}(P,Q)}{\partial n_Q}d\nu(Q);$$
$$B^{(k)}_2(P)=\int_{C_n(\Omega; 2^{k-1}, 2^{k+2})}G^a_{\Omega}(P, Q)d\mu(Q)+\int_{S_n(\Omega; 2^{k-1}, 2^{k+2})}\frac{\partial
G^a_{\Omega}(P,Q)}{\partial n_Q}d\nu(Q)$$ and
$$B^{(k)}_3(P)=\int_{C_n(\Omega; 2^{k+2}, \infty)}G^a_{\Omega}(P, Q)d\mu(Q)$$
$$+\int_{S_n(\Omega; 2^{k+2}, \infty)}\frac{\partial G^a_{\Omega}(P,Q)}{\partial n_Q}d\nu(Q)~~(P\in C_n(\Omega); k=1, 2, 3,\cdots).$$
First we shall prove the existence of an integer $N$ such that
$$H_{\upsilon}\cap I_k(\Omega)\subset\{P=(r, \Theta)\in I_k(\Omega); B^{(k)}_2(P)\geq\frac{1}{2}V(r)\}\eqno{(4.17)}$$
for any integer $k(k\geq N)$. Since $\upsilon(P)$ is finite almost
everywhere on $C_n(\Omega)$, we may apply Lemmas 3 and 4 to
$$\int_{C_n(\Omega)}G^a_{\Omega}(P, Q)d\mu(Q)~and~\int_{S_n(\Omega)}\frac{\partial G^a_{\Omega}(P,Q)}{\partial n_Q}d\nu(Q),$$
respectively, then we can take an integer $N$ such that for any
integer $k(k\geq N)$
$$\frac{W(2^{k-1})}{V(2^{k-1})}\int_{C_n(\Omega; 0, 2^{k-1})}V(t)\varphi(\Phi)d\mu(t, \Phi)\leq\frac{1}{12J_{\Omega}C_2};\eqno{(4.18)}$$
$$\int_{C_n(\Omega; 2^{k+2}, \infty)}W(t)\varphi(\Phi)d\mu(t, \Phi)\leq\frac{1}{12J_{\Omega}C_2};\eqno{(4.19)}$$
$$\frac{W(2^{k-1})}{V(2^{k-1})}\int_{S_n(\Omega; 0, 2^{k-1})}V(t)t^{-1}\frac{\partial\varphi(\Phi)}{\partial n_{\Phi}}d\nu(t, \Phi)\leq\frac{1}{12J_{\Omega}C_2}\eqno{(4.20)}$$
and
$$\int_{S_n(\Omega; 2^{k+2}, \infty)}W(t)t^{-1}\frac{\partial\varphi(\Phi)}{\partial n_{\Phi}}d\nu(t, \Phi)\leq\frac{1}{12J_{\Omega}C_2},\eqno{(4.21)}$$
where
$$J_{\Omega}=\sup_{\Theta\in\Omega}\varphi(\Theta).\eqno{(4.22)}$$
Then for any $P=(r, \Theta)\in I_k(\Omega)(k\geq N)$, we have
$$B^{(k)}_1(P)\leq C_2J_{\Omega}W(r)\int_{C_n(\Omega; 0, 2^{k-1})}V(t)\varphi(\Phi)d\mu(t, \Phi)$$
$$+C_2J_{\Omega}W(r)\int_{S_n(\Omega; 0, 2^{k-1})}V(t)t^{-1}\frac{\partial\varphi(\Phi)}{\partial n_{\Phi}}d\nu(t, \Phi)\leq\frac{V(r)}{6}$$
from (2.1 or 2.2), (3.1 or 3.2), (4.18) and (4.20), and
$$B^{(k)}_3(P)\leq C_2J_{\Omega}V(r)\int_{C_n(\Omega; 2^{k+2}, \infty)}W(t)\varphi(\Phi)d\mu(t, \Phi)$$
$$+C_2J_{\Omega}V(r)\int_{S_n(\Omega; 2^{k+2}, \infty)}W(t)t^{-1}\frac{\partial\varphi(\Phi)}{\partial n_{\Phi}}d\nu(t, \Phi)\leq\frac{V(r)}{6}$$
from (2.1 or 2.2), (3.1 or 3.2), (4.19) and (4.21). Further we can
assume that
$$6\kappa c_o(\upsilon, a)J_{\Omega}\leq V(r)W^{-1}(r)$$
for any $P=(r, \Theta)\in I_k(\Omega)(k\geq N)$. Hence if $P=(r,
\Theta)\in I_k(\Omega)\cap H_{\upsilon}(k\geq N)$, we obtain
$$B^{(k)}_2(P)\geq \upsilon(P)-\frac{V(r)}{6}-B^{(k)}_1(P)-B^{(k)}_3(P)\geq\frac{V(r)}{2}$$
from (4.16) which gives (4.17).

Next we observe from (4.15) and (4.17) that
$$B^{(k)}_2(P)\geq \frac{1}{2}V(2^k)~~(k\geq N)$$
for any $P\in E_k$. If we define a function $u_k(P)$ on
$C_n(\Omega)$ by
$$u_k(P)=2V^{-1}(2^k)B^{(k)}_2(P),$$
then
$$u_k(P)\geq1~~(P\in E_k, k\geq N)$$
and
$$u_k(P)=\int_{C_n(\Omega)}G^a_{\Omega}(P, Q)d\mu_k(Q)+\int_{S_n(\Omega)}\frac{\partial
G^a_{\Omega}(P,Q)}{\partial n_Q}d\nu_k(Q)$$ with two measures
$$
d\mu_k(Q)=\left\{\begin{array}{ll}
 2V^{-1}(2^{k})d\mu(Q) & ~~~~~~~~   \ \ \ (Q\in C_n(\Omega; (2^{k-1}, 2^{k+2}))),  \\
 0

  &~~~~~~~~
 \ \ \ (Q\in C_n(\Omega; (0, 2^{k-1}))\cup C_n(\Omega; [2^{k+2}, \infty)))
 \end{array}\right.
$$
and
$$
d\nu_k(Q)=\left\{\begin{array}{ll}
 2V^{-1}(2^{k})d\nu(Q) & ~~~~~~~~   \ \ \ (Q\in S_n(\Omega; (2^{k-1}, 2^{k+2}))),  \\
 0
  &~~~~~~~~
 \ \ \ (Q\in S_n(\Omega; (0, 2^{k-1}))\cup S_n(\Omega; [2^{k+2},
 \infty))).
 \end{array}\right.
$$
Hence by applying Lemma 7 to $u_k(P)$, we obtain
$$\lambda^a_{\Omega}(E_k)\leq 2V^{-1}(2^{k})\int_{C_n(\Omega; 2^{k-1}, 2^{k+2})}V(t)\varphi(\Phi)d\mu(t, \Phi)$$
$$+2V^{-1}(2^{k})\int_{S_n(\Omega; 2^{k-1}, 2^{k+2})}V(t)t^{-1}\frac{\partial\varphi(\Phi)}{\partial n_{\Phi}}d\nu(t, \Phi)~~(k\geq N).$$
Finally we have by (1.7), (1.8) and (1.9)
$$\sum^{\infty}_{k=N}W(2^k)\lambda^a_{\Omega}(E_k)\lesssim\int_{C_n(\Omega; 2^{N-1}, \infty)}W(t)\varphi(\Phi)d\mu(t, \Phi)$$
$$+\int_{S_n(\Omega; 2^{N-1}, \infty)}W(t)t^{-1}\frac{\partial\varphi(\Phi)}{\partial n_{\Phi}}d\nu(t, \Phi).$$
If we take a sufficiently large $N$, then the integrals of the right
side are finite from Lemmas 3 and 4.

We suppose that a subset $E$ of $C_n(\Omega)$ satisfies
$$\sum^{\infty}_{k=0}W(2^k)\lambda^a_{\Omega}(E_k)<\infty,$$
then we apply the second part of Lemma 7  to $E_k$ and get
$$\sum^{\infty}_{k=1}W(2^k)\left\{\int_{C_n(\Omega)}V(t)\varphi(\Phi)d\mu^*_k(t, \Phi)+\int_{S_n(\Omega)}V(t)t^{-1}
\frac{\partial\varphi(\Phi)}{\partial n_{\Phi}}d\nu^*_k(t,
\Phi)\right\}<\infty,\eqno{(4.23)}$$ where $\mu^*_k$ and $\nu^*_k$
are two positive measures on $C_n(\Omega)$ and $S_n(\Omega)$,
respectively, such that
$$\widehat{R}^{E_k}_1(P)=\int_{C_n(\Omega)}G^a_{\Omega}(P, Q)d\mu^*_k(Q)+\int_{S_n(\Omega)}\frac{\partial
G^a_{\Omega}(P,Q)}{\partial n_Q}d\nu^*_k(Q).\eqno{(4.24)}$$ Consider
a function $\upsilon_0(P)$ on $C_n(\Omega)$ defined by
$$\upsilon_0(P)=\sum^{\infty}_{k=-1}V(2^{k+1})\widehat{R}^{E_k}_1(P)~~(P\in C_n(\Omega)),$$
where
$$E_{-1}=E\cap\{P=(r, \Theta)\in C_n(\Omega); 0<r<1\},$$
then $\upsilon_0(P)$ is a superfunction or identically $\infty$ on
$C_n(\Omega)$. We take any positive integer $k_0$ and represent
$\upsilon_0(P)$ by
$$\upsilon_0(P)=\upsilon_1(P)+\upsilon_2(P),$$
where
$$\upsilon_1(P)=\sum^{k_0+1}_{k=-1}V(2^{k+1})\widehat{R}^{E_k}_1(P),~~~\upsilon_2(P)=
\sum^{\infty}_{k=k_0+2}V(2^{k+1})\widehat{R}^{E_k}_1(P).$$ Since
$\mu^*_k$ and $\nu^*_k$ are concentrated on $B_{E_k}\subset
\overline{E_k}\cap C_n(\Omega)$ and $B'_{E_k}\subset
\overline{E_k}\cap S_n(\Omega)$, respectively, we have from (2.1 or
2.2), (3.1 or 3.2), (1.7) and (1.8) that
$$\int_{C_n(\Omega)}G^a_{\Omega}(P', Q)d\mu^*_k(Q)\leq C_2V(r')\varphi(\Theta')\int_{C_n(\Omega)}W(t)\varphi(\Phi)d\mu^*_k(t, \Phi)$$
$$\leq C_2W(2^{k})V^{-1}(2^{k})V(r')\varphi(\Theta')\int_{C_n(\Omega)}V(t)\varphi(\Phi)d\mu^*_k(t, \Phi)$$
and
$$\int_{S_n(\Omega)}\frac{\partial G^a_{\Omega}(P',Q)}{\partial n_Q}d\nu^*_k(Q)\leq C_2W(2^{k})V^{-1}(2^{k})V(r')\varphi(\Theta')\int_{S_n(\Omega)}V(t)t^{-1}
\frac{\partial\varphi(\Phi)}{\partial n_{\Phi}}d\nu^*_k(t, \Phi)$$
for a point $P'=(r', \Theta')\in C_n(\Omega)(r'\leq 2^{k_0+1}$ and
any integer $k\leq k_0+2)$. Hence we know by (1.7), (1.8) and (1.9)
$$\upsilon_2(P')\lesssim V(r')\varphi(\Theta')\sum^{\infty}_{k=k_0+2}W(2^{k})\int_{C_n(\Omega)}V(t)\varphi(\Phi)d\mu^*_k(t, \Phi)$$
$$+V(r')\varphi(\Theta')\sum^{\infty}_{k=k_0+2}W(2^{k})\int_{S_n(\Omega)}V(t)t^{-1} \frac{\partial\varphi(\Phi)}{\partial n_{\Phi}}d\nu^*_k(t, \Phi).\eqno{(4.25)}$$
This and (4.23) show that $\upsilon_2(P')$ is finite and hence
$\upsilon_0(P)$ is a positive superfunction on $C_n(\Omega)$. To see
$$c(\upsilon_0, a)=\inf_{P\in C_n(\Omega)}\frac{\upsilon_0(P)}{M^a_{\Omega}(P, \infty)}=0,\eqno{(4.26)}$$
we consider the representations of $\upsilon_0(P)$, $\upsilon_1(P)$
and $\upsilon_2(P)$ by Lemma 6
$$\upsilon_0(P)=c(\upsilon_0, a)M^a_{\Omega}(P, \infty)+c_o(\upsilon_0, a)M^a_{\Omega}(P, O)+\int_{C_n(\Omega)}G^a_{\Omega}(P,
Q)d\mu_{(0)}(Q)$$
$$+\int_{S_n(\Omega)}\frac{\partial G^a_{\Omega}(P,Q)}{\partial n_Q}d\nu_{(0)}(Q);$$
$$\upsilon_1(P)=c(\upsilon_1, a)M^a_{\Omega}(P, \infty)+c_o(\upsilon_1, a)M^a_{\Omega}(P, O)+\int_{C_n(\Omega)}G^a_{\Omega}(P,
Q)d\mu_{(1)}(Q)$$
$$+\int_{S_n(\Omega)}\frac{\partial G^a_{\Omega}(P,Q)}{\partial n_Q}d\nu_{(1)}(Q)$$ and
$$\upsilon_2(P)=c(\upsilon_2, a)M^a_{\Omega}(P, \infty)+c_o(\upsilon_2, a)M^a_{\Omega}(P, O)+\int_{C_n(\Omega)}G^a_{\Omega}(P,
Q)d\mu_{(2)}(Q)$$
$$+\int_{S_n(\Omega)}\frac{\partial G^a_{\Omega}(P,Q)}{\partial n_Q}d\nu_{(2)}(Q).$$ It is evident from
(4.24)that $c(\upsilon_1, a)=0$ for any $k_0$. Since $c(\upsilon_0,
a)=c(\upsilon_2, a)$ and
$$c(\upsilon_2, a)=\inf_{P\in C_n(\Omega)}\frac{\upsilon_2(P)}{M^a_{\Omega}(P, \infty)}\leq
\frac{\upsilon_2(P')}{M^a_{\Omega}(P', \infty)}\lesssim
\sum^{\infty}_{k=k_0+2}W(2^{k})\int_{C_n(\Omega)}V(t)\varphi(\Phi)d\mu^*_k(t,\Phi)$$
$$+\sum^{\infty}_{k=k_0+2}W(2^{k})\int_{S_n(\Omega)}V(t)t^{-1} \frac{\partial\varphi(\Phi)}{\partial n_{\Phi}}d\nu^*_k(t, \Phi)
\rightarrow0~~(k_0\rightarrow\infty)$$ from (4.23) and (4.25), we
know $c(\upsilon_0, a)=0$ which is (4.26). Since
$\widehat{R}^{E_k}_1= 1$ on $B_{E_k}\subset \overline{E_k}\cap
C_n(\Omega)$, we know
$$\upsilon_0(P)\geq V(2^{k+1})\geq V(r)$$
for any $P=(r, \Theta)\in B_{E_k}(k=-1, 0, 1, 2,\cdots)$. If we set
$E'=\cup^{\infty}_{k=-1}B_{E_k}$, then
$$E'\subset H_{\upsilon_0}.\eqno{(4.27)}$$
Since $E'$ is equal to $E$ except a polar set $S$, we can take
another positive superfunction $\upsilon_3$ on $C_n(\Omega)$ such
that $\upsilon_3=G^a_{\Omega}\eta$ with a positive measure $\eta$ on
$C_n(\Omega)$ and $\upsilon_3$ is identically $\infty$ on $S$.
Finally we define a positive superfunction $\upsilon$ on
$C_n(\Omega)$ by
$$\upsilon=\upsilon_0+\upsilon_3.$$
Since $c(\upsilon_3, a)=0$, it is easy to see from (4.26) that
$c(\upsilon, a)=0$. In addition we know from (4.27) that $E\subset
H_{\upsilon}.$ Thus we complete the proof that subset $E$ of
$C_n(\Omega)$ is a-rarefied at $\infty$ with respect to
$C_n(\Omega)$.

\vspace{0.3mm}{\bf Proof of the Theorem 7.} By Lemma 6 we have
$$\upsilon(P)=c(\upsilon, a)M^a_{\Omega}(P, \infty)+c_o(\upsilon, a)M^a_{\Omega}(P, O)
+\int_{C_n(\Omega)}G^a_{\Omega}(P, Q)d\mu(Q)$$
$$+\int_{S_n(\Omega)}\frac{\partial G^a_{\Omega}(P, Q)}{\partial n_Q}d\nu(Q)$$
for a unique positive measure $\mu$ on $C_n(\Omega)$ and a unique
positive measure $\nu$ on $S_n(\Omega)$, respectively. Then
$$\upsilon_1(P)=\upsilon(P)-c(\upsilon, a)M^a_{\Omega}(P, \infty)-c_o(\upsilon, a)M^a_{\Omega}(P, O)
~~(P=(r, \Theta)\in C_n(\Omega))$$ also is a positive superfunction
on $C_n(\Omega)$ such that
$$\inf_{P=(r, \Theta)\in C_n(\Omega)}\frac{\upsilon_1(P)}{M^a_{\Omega}(P, \infty)}=0.$$
We shall prove the existence of a-rarefied set $E$ at $\infty$ with
respect to $C_n(\Omega)$ such that
$$\upsilon_1(P)V^{-1}(r)~~(P=(r, \Theta)\in C_n(\Omega))$$
uniformly converges to 0 on $C_n(\Omega)\diagdown E$ as
$r\rightarrow\infty$. Let $\{\varepsilon_i\}$ be a sequence of
positive numbers $\varepsilon_i$ satisfying
$\varepsilon_i\rightarrow0$ $i\rightarrow\infty$. Put
$$E_i=\{P=(r, \Theta)\in C_n(\Omega); \upsilon_1(P)\geq\varepsilon_iV(r)\}~~(k=1, 2, 3,\cdots).$$
Then $E_i(k=1, 2, 3,\cdots)$ is a-rarefied sets at $\infty$ with
respect to $C_n(\Omega)$ and hence by Theorem 6
$$\sum^{\infty}_{k=0}W(2^k)\lambda^a_{\Omega}((E_i)_k)<\infty~~(i=1, 2, 3,\cdots).$$
We take a sequence $\{q_i\}$ such that
$$\sum^{\infty}_{k=q_i}W(2^k)\lambda^a_{\Omega}((E_i)_k)<\frac{1}{2^i}~~(i=1, 2, 3,\cdots)$$
and set
$$E=\cup^{\infty}_{i=1}\cup^{\infty}_{k=q_i}(E_i)_k.$$
Because $\lambda^a_{\Omega}$ is a countably sub-additive set
function as in H.Aikawa [1] and M.Ess\'{e}n and H.L.Jackson [9],
then
$$\lambda^a_{\Omega}(E_m)\leq\sum^{\infty}_{i=1}\sum^{\infty}_{k=q_i}\lambda^a_{\Omega}(E_i\cap I_k\cap I_m)~~(m=1, 2, 3,\cdots).$$
Since
$$\sum^{\infty}_{m=1}\lambda^a_{\Omega}(E_m)W(2^m)\leq\sum^{\infty}_{i=1}\sum^{\infty}_{k=q_i}\sum^{\infty}_{m=1}\lambda^a_{\Omega}(E_i\cap I_k\cap I_m)W(2^m)$$
$$=\sum^{\infty}_{i=1}\sum_{k=q_i}\lambda^a_{\Omega}((E_i)_k)W(2^k)\leq\sum^{\infty}_{i=1}\frac{1}{2^i}=1,$$
we know by Theorem 6 that $E$ is a a-rarefied set at $\infty$ with
respect to $C_n(\Omega)$. It is easy to see that
$$\upsilon(P)V^{-1}(r)~~(P=(r, \Theta)\in C_n(\Omega))$$ uniformly converges to 0 on $C_n(\Omega)\diagdown E$ as
$r\rightarrow\infty$.

\vspace{0.3mm}{\bf Proof of the Theorem 8.} Since $\lambda_{E_k}$ is
concentrated on $B_{E_k}\subset \overline{E_k}\cap C_n(\Omega)$, we
see
$$\gamma^a_{\Omega}(E_k)=\int_{C_n(\Omega)}\widehat{R}^{E_k}_{M^a_{\Omega}(., \infty)}(P)d\lambda_{E_k}(P)
\leq\int_{C_n(\Omega)}M^a_{\Omega}(P, \infty)d\lambda_{E_k}(P)\leq
J_{\Omega}V(2^{k+1})\lambda^a_{\Omega}(E_k)$$ and hence
$$\sum^{\infty}_{k=0}\gamma^a_{\Omega}(E_k) W(2^k)V^{-1}(2^k)\lesssim\sum^{\infty}_{k=0}W(2^{k})\lambda^a_{\Omega}(E_k),$$
which gives the conclusion in the first part from Theorems 5 and 6.
To prove the second part, we put
$J'_{\Omega}=\min_{\Theta\in\overline{\Omega'}}\varphi(\Theta)$.
Since
$$M^a_{\Omega}(., \infty)=V(r)\varphi(\Theta)\geq J'_{\Omega}V(r)\geq J'_{\Omega}V(2^k)~~(P=(r, \Theta)\in E_k)$$
and
$$\widehat{R}^{E_k}_{M^a_{\Omega}(., \infty)}(P)=M^a_{\Omega}(., \infty)$$
for any $P=(r, \Theta)\in B_{E_k}$, we have
$$\gamma^a_{\Omega}(E_k)=\int_{C_n(\Omega)}\widehat{R}^{E_k}_{M^a_{\Omega}(., \infty)}(P)d\lambda_{E_k}(P)\geq J'_{\Omega}V(2^k)\lambda^a_{\Omega}(E_k).$$
Since
$$J'_{\Omega}\sum^{\infty}_{k=0}\lambda^a_{\Omega}(E_k) W(2^k)\leq\sum^{\infty}_{k=0}V^{-1}(2^k)W(2^{k})\gamma^a_{\Omega}(E_k)<\infty$$
from Theorem 5, it follows from Theorem 6 that $E$ is a-rarefied at
$\infty$ with respect to $C_n(\Omega)$.

\end{document}